# CONVERGENCE OF A FAMILY OF SERIES

Florentin Smarandache, Ph D
Full Professor
Chair of Department of Math & Sciences
University of New Mexico
200 College Road
Gallup, NM 87301, USA
E-mail: smarand@unm.edu

In this article we will construct a family of expressions $\mathcal{E}(n)$. For each element $E(n)$ from $\mathcal{E}(n)$, the convergence of the series $\sum_{n \geq n_E} E(n)$ can be determined in accordance to the theorems of this article.

This article gives also applications.

**(1) Preliminary**

To render easier the expression, we will use the recursive functions. We will introduce some notations and notions to simplify and reduce the size of this article.

**(2) Definitions: lemmas.**

We will construct recursively a family of expressions $\mathcal{E}(n)$.
For each expression $E(n) \in \mathcal{E}(n)$, the degree of the expression is defined recursively and is denoted $d^0 E(n)$, and its dominant (leading) coefficient is denoted $c(E(n))$.

1. If $a$ is a real constant, then $a \in \mathcal{E}(n)$.

   $d^0 a = 0$ and $c(a) = a$.

2. The positive integer $n \in \mathcal{E}(n)$.

   $d^0 n = 1$ and $c(n) = 1$.

3. If $E_1(n)$ and $E_2(n)$ belong to $\mathcal{E}(n)$ with $d^0 E_1(n) = r_1$ and $d^0 E_2(n) = r_2$, $c(E_1(n)) = a_1$ and $c(E_2(n)) = a_2$, then:

   a) $E_1(n) E_2(n) \in \mathcal{E}(n)$; $d^0(E_1(n) E_2(n)) = r_1 + r_2$; $c(E_1(n) E_2(n))$ which is $a_1 a_2$.

   b) If $E_2(n) \neq 0 \ \forall n \in \mathbb{N}(n \geq n_{E_2})$, then $\dfrac{E_1(n)}{E_2(n)} \in \mathcal{E}(n)$ and

$$d^0\left(\frac{E_1(n)}{E_2(n)}\right) = r_1 - r_2, \quad c\left(\frac{E_1(n)}{E_2(n)}\right) = \frac{a_1}{a_2}.$$

c) If $\alpha$ is a real constant and if the operation used is well defined, $(E_1(n))^\alpha$ (for all $n \in N$, $n \geq n_{E_1}$), then:

$$(E_1(n))^\alpha \in \mathcal{E}(n), \quad d^0\left((E_1(n))^\alpha\right) = r_1\alpha, \quad c\left((E_1(n))^\alpha\right) = a_1^\alpha$$

d) If $r_1 \neq r_2$, then $E_1(n) \pm E_2(n) \in \mathcal{E}(n)$, $d^0(E_1(n) \pm E_2(n))$ is the max of $r_1$ and $r_2$, and $c(E_1(n) \pm E_2(n)) = a_1$, respectively $a_2$ resulting that the grade is $r_1$ and $r_2$.

e) If $r_1 = r_2$ and $a_1 + a_2 \neq 0$, then $E_1(n) + E_2(n) \in \mathcal{E}(n)$, $d^0(E_1(n) + E_2(n)) = r_1$ and $c(E_1(n) + E_2(n)) = a_1 + a_2$.

f) If $r_1 = r_2$ and $a_1 - a_2 \neq 0$, then $E_1(n) - E_2(n) \in \mathcal{E}(n)$, $d^0(E_1(n) - E_2(n)) = r_1$ and $c(E_1(n) - E_2(n)) = a_1 - a_2$.

4. All expressions obtained by applying a finite number of step 3 belong to $\mathcal{E}(n)$.

**Note 1.** From the definition of $\mathcal{E}(n)$ it results that, if $E(n) \in \mathcal{E}(n)$ then $c(E(n)) \neq 0$, and that $c(E(n)) = 0$ if and only if $E(n) = 0$.

**Lemma 1.** If $E(n) \in \mathcal{E}(n)$ and $c(E(n)) > 0$, then there exists $n' \in \mathbb{N}$, such that for all $n > n'$, $E(n) > 0$.

*Proof:* Let's consider $c(E(n)) = a_1 > 0$ and $d^0(E(n)) = r$.

If $r > 0$, then $\lim_{n \to \infty} E(n) = \lim_{n \to \infty} n^r \frac{E(n)}{n^r} = \lim_{n \to \infty} a_1 n^r = +\infty$, thus there exists $n' \in \mathbb{N}$ such that, for $n > n'$ we have $E(n) > 0$.

If $r < 0$, then $\lim_{n \to \infty} \frac{1}{E(n)} = \lim_{n \to \infty} \frac{n^{-r}}{\frac{E(n)}{n^r}} = \frac{1}{a_1} \lim_{n \to \infty} n^{-r} = +\infty$ thus there exists $n' \in \mathbb{N}$, such that for all $n > n'$, $\frac{1}{E(n)} > 0$, hence we have $E(n) > 0$.

If $r = 0$, then $E(n)$ is a positive real constant, or $\frac{E_1(n)}{E_2(n)} = E(n)$, with $d^0 E_1(n) = d^0 E_2(n) = r_1 \neq 0$, according to what we have just seen, $c\left(\frac{E_1(n)}{E_2(n)}\right) = \frac{c(E_1(n))}{c(E_2(n))} = c(E(n)) > 0$.

Then: $c(E_1(n)) > 0$ and $c(E_2(n)) < 0$: it results

there exists $n_{E_1} \in \mathbb{N}$, $\forall n \in \mathbb{N}$ and $n \geq n_{E_1}$, $E_1(n) > 0$

there exists $n_{E_2} \in \mathbb{N}$, $\forall n \in \mathbb{N}$ and $n \geq n_{E_2}$, $E_2(n) > 0$ $\Rightarrow$

there exists $n_E = \max(n_{E_1}, n_{E_2}) \in \mathbb{N}$, $\forall n \in \mathbb{N}$, $n \geq n_E$, $E(n) \dfrac{E_1(n)}{E_2(n)} > 0$

then $c(E_1(n)) < 0$ and $c(E_2(n)) < 0$ and it results:

$$E(n) = \frac{E_1(n)}{E_2(n)} = \frac{-E_1(n)}{-E_2(n)}$$ which brings us back to the precedent case.

**Lemma 2:** If $E(n) \in \mathcal{E}(n)$ and if $c(E(n)) < 0$, then it exists $n' \in \mathbb{N}$, such that qqst $n > n'$, $E(n) < 0$.

*Proof:*

The expression $-E(n)$ has the propriety that $c(-E(n)) > 0$, according to the recursive definition. According to lemma 1: there exists $n' \in \mathbb{N}$, $n \geq n'$, $-E(n) > 0$, i.e. $+E(n) < 0$, q. e. d.

**Note 2.** To prove the following theorem, we suppose known the criterion of convergence of the series and certain of its properties

### (3) Theorem of convergence and applications.

**Theorem:** Let's consider $E(n) \in \mathcal{E}(n)$ with $d^0(E(n)) = r$ having the series

$$\sum_{n \geq n_\varepsilon} E(n), \quad E(n) \not\equiv 0.$$

Then:

A) If $r < -1$ the series is absolutely convergent.
B) If $r \geq -1$ it is divergent where $E(n)$ is well defined $\forall n \geq n_E, n \in \mathbb{N}$.

*Proof:* According to lemmas 1 and 2, and because:

the series $\sum_{n \geq n_E} E(n)$ converge $\Leftrightarrow$ the series $-\sum_{n \geq n_E} E(n)$ converge,

we can consider the series $\sum_{n \geq n_E} E(n)$ like a series with positive terms.

We will prove that the series $\sum_{n \geq n_E} E(n)$ has the same nature as the series $\sum_{n \geq 1} \dfrac{1}{n^{-r}}$.

Let us apply the second criterion of comparison:

$$\lim_{n \to \infty} \frac{E(n)}{\dfrac{1}{n^{-r}}} = \lim_{n \to \infty} \frac{E(n)}{n^r} = c(E(n)) \neq \pm \infty.$$

According to the note 1 if $E(n) \not\equiv 0$ then $c(E(n)) \neq 0$ and then the series $\sum_{n \geq n_E} E(n)$ has the same nature as the series $\sum_{n \geq 1} \dfrac{1}{n^{-r}}$, i.e.:

A) If $r < -1$ then the series is convergent;
B) If $r > -1$ then the series is divergent;

For $r < -1$ the series is absolute convergent because it is a series with positive terms.

**Applications:**

We can find many applications of these. Here is an interesting one:
If $P_q(n)$, $R_s(n)$ are polynomials in $n$ of degrees $q, s$ respectively, and that $P_q(n)$ and $R_s(n)$ belong to $\mathcal{E}(n)$:

1) $\sum_{n \geq n_{PR}} \dfrac{\sqrt[k]{P_q(n)}}{\sqrt[h]{R_s(n)}}$ is $\begin{cases} \text{convergent, if } s/h - q/k > 1 \\ \text{divergent, if } s/h - q/k \leq 1 \end{cases}$

2) $\sum_{n \geq n_R} \dfrac{1}{R_s(n)}$ is $\begin{cases} \text{convergent, if } s > 1 \\ \text{divergent, if } s \leq 1 \end{cases}$

**Example:** The series $\sum_{n \geq 2} \dfrac{\sqrt[2]{n+1} \cdot \sqrt[3]{n-7} + 2}{\sqrt[5]{n^2} - 17}$ is divergent because $\dfrac{2}{5} - \left(\dfrac{1}{2} + \dfrac{1}{3}\right) < 1$ and if we call $E(n)$ the quotient of this series, $E(n)$ belongs to $\mathcal{E}(n)$ and it is well defined for $n \geq 2$.